\documentclass[smallextended]{svjour3}
\smartqed
\usepackage{graphicx}
\usepackage{amsmath}
\newtheorem{Theorem}{Theorem}
\newtheorem{Proposition}{Proposition}
\newtheorem{Definition}{Definition}
\newtheorem{Corollary}{Corollary}
\newtheorem{Lemma}{Lemma}

\begin{document}
\title{$3$-Kenmotsu Manifolds}
%\subtitle{$3$-Kenmotsu Manifolds}
\author{Hassan Attarchi}
\institute{H. Attarchi \at
              School of Mathematics \\ Georgia Institute of Technology \\
              \email{hattarchi@gatech.edu}}
\date{Received: date / Accepted: date}
\maketitle
\begin{abstract}
In this paper, a $3$-Kenmotsu structure is defined on a $4n+1$ dimensional manifold where such structure seems to be never studied before.
\end{abstract}
%%%%%%%%%%%%%%%%%%%%%%%%%%%%%%%%%%%%%%%%%%%%%%%%%%%%%%%%%%%%%%%%%%%%%%%%%
\section{Introduction}
A $(2n+1)$-dimensional smooth manifold $M$ is said to have an almost contact structure if the structural group of its tangent bundle reduces to $U(n)\times 1$, where $U(n)$ is the unitary group of degree $n$ \cite{Gray59}. Equivalently, an almost contact structure is given by a triple $(\phi,\xi,\eta)$ satisfying certain conditions \cite{Sasaki60,Sasaki61,Sasaki62}. It is well-known that the almost contact metric structure on an odd-dimensional manifold is analogous to an almost Hermitian structure on an even-dimensional manifold. For example, in \cite{Tashiro63}, it is shown that a hypersurface in an almost Hermitian manifold has an almost contact metric structure. Moreover, a Kenmotsu manifold is a locally warped product of an interval of the real line and a Kahler manifold with a special warped function \cite{Kenmotsu72}. Also, in \cite{Ganchev08}, it is proven that the warped product of a Sasakian manifold with the real line under a certain conformal change will result in a Kahlerian structure.\par
However, the situation is completely different when someone wants to construct an odd-dimensional structure analogous to the quaternion spaces. There are two choices, $(4n+1)$ and $(4n-1)$-dimensional spaces. In literature, what is known as an almost contact 3-structure is based on a $(4n-1)$-dimensional space or equivalently a $(4n+3)$-dimensional space \cite{Kuo70}. In this case, the structural group of almost contact 3-structure reduces to $Sp(n)\times I_3$ \cite{Ishihara74,Kuo70}.\par
On the other hand, what makes this difference more deepening is about constructing an analogous structure of Kenmotsu manifolds on $(4n+3)$-dimensional spaces. The $(4n+3)$-dimensional manifolds as underlying manifolds of almost contact $3$-structures lead us to have only almost contact $3$-structure or $3$-Sasakian manifold \cite{Blair02,Kashiwada01}. This terminology is completely in agreement with the result of \cite{Ganchev08} which was mentioned earlier. This means one can expect that the warped product of a $3$-Sasakian manifold with the real line under a certain conformal change will result in a quaternion Kahler manifold. While it cannot support the strategy behind constructing the Kenmotsu manifolds \cite{Kenmotsu72}. To have a structure which is locally the warped product of an interval and a quaternion Kahler manifold, one should study a new structure on $(4n+1)$-dimensional manifolds where its structural group reduces to $Sp(n)\times I$. To the best of my knowledge, nobody has considered and studied $(4n+1)$-dimensional manifolds from this perspective. To name one of the most recent works on $(4n+1)$-dimensional manifolds, one can see \cite{Galaev18}.\par
The structure of the paper is the following. In Section 2, we introduce necessary notations and study some properties of Kenmotsu manifolds. Section 3 deals with the definition of a $3$-Kenmotsu manifold as an odd-dimensional analogous structure of the quaternion Kahler manifolds. In Section 4, it is proved that $3$-Kenmotsu manifolds are Einstein spaces. Moreover, it is shown that the $\varphi_\alpha-$holomorphic sectional curvature $\mathbf{H}_\alpha(X)$ of this structure satisfies the equation $\sum_{\alpha=1}^3\mathbf{H}_\alpha(X)=-3$ for all $X\in\Gamma H$ where $H$ is the contact distribution of Kenmotsu structures. In Section 5, an example of $3$-Kenmotsu manifolds is studied.
%%%%%%%%%%%%%%%%%%%%%%%%%%%%%%%%%%%%%%%%%%%%%%%%%%%%%%%%%%%%%%%%%%%%%%%%%
\section{Preliminaries and Notations}
Let $(\varphi,\eta,\xi,g)$ be an almost contact metric structure on $(2n+1)$-dimensional manifold $M$, where $\varphi\in End(TM)$, $\xi$ is the Reeb vector field, and $\eta$ is its dual 1-form with respect to the Riemannian metric $g$. Also, they satisfy following properties,
\begin{equation}\label{0}
\begin{split}
\varphi^2X=-X&+\eta(X)\xi,\ \ \eta(\xi)=1,\ \ \eta(\varphi)=0,\ \ \varphi\xi=0,\\
& g(\varphi X,\varphi Y)=g(X,Y)-\eta(X)\eta(Y),
\end{split}
\end{equation}
where $X,Y\in\Gamma TM$ \cite{Blair02}. An almost contact metric manifold is called a Kenmotsu manifold if
\begin{equation}\label{1}
(\nabla_{X}\varphi)Y=g(\varphi X,Y)\xi-\eta(Y)\varphi X,
\end{equation}
where $\nabla$ is the Levi-Civita connection of the Riemannian metric $g$ and $X,Y\in\Gamma TM$. This implies that:
\begin{equation}\label{2}
\nabla_{X}\xi=X-\eta(X)\xi,\ \ \ \ \ (\nabla_{X}\eta)Y=g(X,Y)-\eta(X)\eta(Y),
\end{equation}
where $X,Y\in\Gamma TM$ \cite{Kenmotsu72}.
Moreover, the $1$-form $\eta$ of almost contact structure of a Kenmotsu manifold is closed (i.e. $d\eta$=0). The Kahler form $\Omega$ is defined on an almost contact metric manifold as follows:
\begin{equation}\label{omega}
\Omega(X,Y):=g(X,\varphi Y),
\end{equation}
where $X,Y\in\Gamma TM$. Let $(M,\varphi,\eta,\xi,g)$ be a Kenmotsu manifold of dimension $2n+1$. The Kahler form $\Omega$ of $M$ satisfies the equation \cite{Olszak91},
$$d\Omega=\eta\wedge\Omega.$$
Consider the foliation $F$ of the Reeb vector field $\xi$ on $M$. Then, there are local frames $\{U; x^0, x^i\}$ adapted to this foliation where $\xi=\partial/\partial x^0$ on $U$ \cite{Molino88}. The local vector fields of this coordinate system define local frames $\{\partial/\partial x^0, \delta/\delta x^i\}$ on $U$ where,
\begin{equation}\label{4.0}
\frac{\delta}{\delta x^{i}}=\frac{\partial}{\partial x^{i}}-\eta_{i}\frac{\partial}{\partial x^{0}},
\end{equation}
where $\eta_{i}:=\eta(\partial/\partial x^{i})$ for $i=1,...,2n$. The Riemannian metric $g$ in the local frames (\ref{4.0}) has the following format,
\begin{equation}\label{metric2}
g:=\left(\begin{array}{cccc}
1&0\\
0&(g_{ij})
\end{array}\right),
\end{equation}
where $g_{ij}=g(\frac{\delta}{\delta x^{i}},\frac{\delta}{\delta x^{j}})$. It is easy to check that the Lie bracket of these local frames satisfies the following properties,
\begin{equation}\label{4.3}
[\frac{\partial}{\partial x^{0}},\frac{\delta}{\delta x^{i}}]=0,\ \ \ \ [\frac{\delta}{\delta x^{i}},\frac{\delta}{\delta x^{j}}]=0.
\end{equation}
\begin{Lemma}\label{xig}
	Let $(M,\varphi,\eta,\xi,g)$ be a Kenmotsu manifold, then components of $(g_{ij})$ defined in (\ref{metric2}) satisfy the following equation:
	$$\xi g_{ij}=2g_{ij}.$$
\end{Lemma}
\proof
	Using (\ref{2}), (\ref{4.3}) and the Levi-Civita connection $\nabla$ on $M$, they imply that,
	\begin{equation*}
	\begin{split}
	\xi g_{ij}=&\xi g(\frac{\delta}{\delta x^i},\frac{\delta}{\delta x^j})=g(\nabla_\xi\frac{\delta}{\delta x^i},\frac{\delta}{\delta x^j})+g(\frac{\delta}{\delta x^i},\nabla_\xi\frac{\delta}{\delta x^j})\\
	=& g(\nabla_\frac{\delta}{\delta x^i}\xi,\frac{\delta}{\delta x^j})+g([\xi,\frac{\delta}{\delta x^i}],\frac{\delta}{\delta x^j})+g(\frac{\delta}{\delta x^i},\nabla_\frac{\delta}{\delta x^j}\xi)+g(\frac{\delta}{\delta x^i},[\xi,\frac{\delta}{\delta x^j}])\\
	=& g(\frac{\delta}{\delta x^i},\frac{\delta}{\delta x^{j}})+g(\frac{\delta}{\delta x^i},\frac{\delta}{\delta x^{j}})=2g(\frac{\delta}{\delta x^i},\frac{\delta}{\delta x^{j}})=2g_{ij}
	\end{split}
	\end{equation*}\qed
\begin{Proposition}\label{Levi-Civita}
	The Levi-Civita connection $\nabla$ on a Kenmotsu manifold $M$ has the local components:
	\begin{equation*}
	\begin{split}
	\nabla_{\frac{\delta}{\delta x^i}}\frac{\delta}{\delta x^j}=&\Gamma_{ij}^k\frac{\delta}{\delta x^k}-g_{ij}\frac{\partial}{\partial x^0},\\
	\nabla_{\frac{\partial}{\partial x^0}}\frac{\delta}{\delta x^i}=&\nabla_{\frac{\delta}{\delta x^i}}\frac{\partial}{\partial x^0}=\frac{\delta}{\delta x^i},\\
	\nabla_{\frac{\partial}{\partial x^{0}}}\frac{\partial}{\partial x^{0}}=& 0,
	\end{split}
	\end{equation*}
	where $\Gamma_{ij}^k=\frac{1}{2}g^{kl}\{\frac{\delta g_{lj}}{\delta x^i}+\frac{\delta g_{il}}{\delta x^j}-\frac{\delta g_{ij}}{\delta x^l}\}$.
\end{Proposition}
\proof
	The proof follows from (\ref{4.0})--(\ref{4.3}) and Lemma \ref{xig}.\qed
%%%%%%%%%%%%%%%%%%%%%%%%%%%%%%%%%%%%%%%%%%%%%%%%%%%%%%%%%%%%%%%%%%%%%%%%%
\section{$3$-Kenmotsu Manifolds}
A $3$-Kenmotsu manifold is defined as follows:
\begin{Definition}\label{Def1}
	Let $M$ be an $m$-dimensional smooth manifold. Then, $M$ is a $3$-Kenmotsu manifold if it is equipped with three Kenmotsu structures $(\varphi_\alpha,\eta,\xi,g)$, where $\alpha=1,2,3$ and
	\begin{equation}\label{def3ken}
	\varphi_k=\varphi_i\circ\varphi_j,
	\end{equation}
	for all $(i,j,k)$ where they are even permutations of $(1,2,3)$.
\end{Definition}
It is easy to see that any $3$-Kenmotsu manifold is a $(4n+1)$-dimensional manifold. The transversal distribution to the foliation $F$ given by Reeb vector field $\xi$ is integrable because $\eta$ is a closed form. Let denote the transversal distribution by $H$, then:
\begin{equation}\label{defH}
	H=\{X\in TM\ |\ \eta(X)=0\}.
\end{equation}
One can define three almost complex structures $J_\alpha=\varphi_\alpha|_H$ for $\alpha=1,2,3$ on the maximal integral submanifolds of $H$. It will be easy to check that $J_\alpha$ for $\alpha=1,2,3$ form an almost quaternion structure on these maximal integral submanifolds of $H$. Thus, the distribution $H$ is a $4n$-dimensional distribution and consequently any $3$-Kenmotsu manifold is $(4n+1)$ dimensional manifold.\par
Moreover, there is a natural volume form on a $3$-Kenmotsu manifold given by:
$$V\!ol_{_M}=\Omega^n\wedge\eta,$$
where $\Omega$ is a $4$-form of maximum rank with the following structure,
$$\Omega=\Omega_1\wedge\Omega_1+\Omega_2\wedge\Omega_2+\Omega_3\wedge\Omega_3,$$
where $\Omega_\alpha(.,.)=g(.,\varphi_\alpha.)$ for $\alpha=1,2,3$.
\begin{Theorem}\label{thm1}
	Let $M$ be a $3$-Kenmotsu manifold. Then, for any $p\in M$, some neighborhood $U$ of that point is identified with a warped product space $(-\epsilon,+\epsilon)\times_fV$ such that $(-\epsilon,+\epsilon)$ is an open interval, $f(t)=ce^t$, and $V$ is a quaternion Kahler manifold.
\end{Theorem}
\proof
	All three Kenmotsu structures of a $3$-Kenmotsu manifold share the same Reeb vector field, therefore they have the same distribution $H$ defined in (\ref{defH}). Theorem 4 in \cite{Kenmotsu72} implies that for each Kenmotsu structure $(\varphi_\alpha,\eta,\xi,g)$ on $M$, the manifold $M$ is locally warped product of an open interval and a Kahler manifold with the warped function $f(t)=ce^t$. This means a $3$-Kenmotsu manifold is locally warped product of an open interval and an almost quaternion manifold. Considering the Levi-Civita connection $\bar{\nabla}$ on a maximal integral submanifold of the foliation $H$ which is locally the same as the almost quaternion manifold in the warped product. By using Eqs. (\ref{1}) and (\ref{omega}), it will be easy to check that $\bar{\nabla}J_\alpha=0$ and $\bar{\nabla}\Omega_\alpha=0$ for $\alpha=1,2,3$. Thus, theorem 1.1 in \cite{Ishihara74} implies that the structure on the maximal integral submanifolds of the foliation $H$ are quaternion Kahler manifolds.\qed
\begin{Corollary}
	The structural group of the tangent bundle of a $3$-Kenmotsu manifold will be reducible to $Sp(n)\times I$.
\end{Corollary}
Similar to the quaternion structures, one can show that there is no fourth Kenmotsu structure $(\varphi_4,\eta,\xi,g)$ on $M$ which satisfies the anti-commutativity conditions with the other three structures \cite{Blair02,Tachibana70}. To see this, let $J_\alpha=\varphi_\alpha|_H$ for $\alpha=1,2,3,4$ be almost complex structures induced on the maximal integral submanifolds of $H$. Then, $J_i\circ J_4=-J_4\circ J_i$ for $i=1,2,3$ and
$$J_3\circ J_4=J_1\circ J_2\circ J_4=-J_1\circ J_4\circ J_2=J_4\circ J_1\circ J_2=J_4\circ J_3,$$
which is a contradiction with the anti commutativity condition.\par
Also, it is well-known if there are two almost complex structure (or two Sasakian structure) on $M$, then under some simple conditions, one can construct a quaternion structure (or $3$-Sasakian structure) \cite{Blair02,Kuo70} based on those two structures.
\begin{Proposition}\label{Pro1}
	Assume $M$ is a $(4n+1)$-dimensional differential manifold. If there are two Kenmotsu structures $(\varphi_1,\eta,\xi,g)$ and $(\varphi_2,\eta,\xi,g)$ on $M$ satisfying $\varphi_1\circ\varphi_2=-\varphi_2\circ\varphi_1$, then $M$ will have a 3-Kenmotsu structure.
\end{Proposition}
\proof
	Let $\varphi_3=\varphi_1\circ\varphi_2$. Then to prove $(\varphi_3,\eta,\xi,g)$ is a Kenmotsu structure on $M$, one can show that $\varphi_3$ satisfies the equation (\ref{1}). Thus,
	\begin{equation*}
	\begin{split}
	(\nabla_X\varphi_3)Y=&(\nabla_X\varphi_1\circ\varphi_2)Y=\nabla_X(\varphi_1\circ\varphi_2Y)-\varphi_1\circ\varphi_2\nabla_XY\\
	=&(\nabla_X\varphi_1)\varphi_2Y+\varphi_1\nabla_X\varphi_2Y-\varphi_1\circ\varphi_2\nabla_XY\\
	=& g(\varphi_1X,\varphi_2Y)\xi-\eta(\varphi_2Y)\varphi_1X+\varphi_1((\nabla_X\varphi_2)Y+\varphi_2\nabla_XY)-\varphi_1\circ\varphi_2\nabla_XY\\
	=& g(\varphi_1X,\varphi_2Y)\xi+\varphi_1(g(\varphi_2X,Y)\xi-\eta(Y)\varphi_2X)\\
	=& g(\varphi_1\circ\varphi_2X,Y)\xi-\eta(Y)\varphi_1\circ\varphi_2(X)=g(\varphi_3X,Y)\xi-\eta(Y)\varphi_3(X).
	\end{split}
	\end{equation*}
	Moreover, it is easy to check that $\varphi_k=\varphi_i\circ\varphi_j$ satisfies for all even permutations $(i,j,k)$ of $(1,2,3)$.\qed
%%%%%%%%%%%%%%%%%%%%%%%%%%%%%%%%%%%%%%%%%%%%%%%%%%%%%%%%%%%%%%%%%%%%%%%%%
\section{Some Properties of a $3$-Kenmotsu Manifold}
Let $\bar{M}$ be a maximal integral submanifold of the foliation $H$ in the $3$-Kenmotsu manifold $M$. Let denote the Levi-Civita connections of $M$ and $\bar{M}$ by $\nabla$ and $\bar{\nabla}$, respectively. The relation between these two connections is given by the Gauss formula \cite{Lee97} as follows:
\begin{equation}\label{gauss}
\nabla_XY=\bar{\nabla}_XY+h(X,Y),
\end{equation}
where $X,Y\in\Gamma T\bar{M}$. It follows from Proposition \ref{Levi-Civita} that,
\begin{equation}\label{h=g}
h(X,Y)=g(X,Y)\xi,
\end{equation}
where $X,Y\in\Gamma T\bar{M}$. Let $R$ and $\bar{R}$ be the curvature tensors of $M$ and $\bar{M}$, respectively. Then (\ref{gauss}) and (\ref{h=g}) imply that $R$ and $\bar{R}$ satisfy the following equation:
\begin{equation}\label{cur1}
\begin{split}
R(X,Y,Z,W)=&\bar{R}(X,Y,Z,W)-g(h(X,W),h(Y,Z))+g(h(X,Z),h(Y,W))\\
=&\bar{R}(X,Y,Z,W)-g(X,W)g(Y,Z)+g(X,Z)g(Y,W),
\end{split}
\end{equation}
where $X,Y,Z,W\in\Gamma T\bar{M}$. Moreover, Ricci tensors $Ric$ and $\bar{R}ic$ of curvature tensors $R$ and $\bar{R}$, respectively, satisfy the following equation for all $X,Y\in\Gamma T\bar{M}$:
\begin{equation}\label{cur2}
\begin{split}
Ric(X,Y)=&\sum_{k=1}^{4n}R(E_k,X,Y,E_k)+R(\xi,X,Y,\xi)\\ =&\sum_{k=1}^{4n}R(E_k,X,Y,E_k)-g(X,Y),
\end{split}
\end{equation}
where $\{E_1,E_2,...,E_{4n}\}$ is an orthonormal local basis of $\Gamma H$. Then, it follows from (\ref{cur1}) and (\ref{cur2}) that,
\begin{equation}\label{ric}
\begin{split}
\bar{R}ic(X,Y)=&\sum_{k=1}^{4n}\bar{R}(E_k,X,Y,E_k)\\
=&\sum_{k=1}^{4n}R(E_k,X,Y,E_k)+\sum_{k=1}^{4n}g(X,Y)-\sum_{k=1}^{4n}g(E_k,Y)g(X,E_k)\\
=& Ric(X,Y)+4ng(X,Y).
\end{split}
\end{equation}
\begin{Theorem}\label{par}
	Let $M$ be a $3$-Kenmotsu manifold of dimension $\geq9$. Then, its Ricci tensor is parallel.
\end{Theorem}
\proof
	Lemma 3.1 in \cite{Ishihara74} implies that the Ricci tensor $\bar{R}ic$ of the quaternion Kahler manifold $\bar{M}$ is parallel. Then, it follows from (\ref{ric}) that the Ricci tensor $Ric$ of $3$-Kenmotsu manifold $M$ is parallel.\qed
\begin{Theorem}\label{Ein}
	Let $M$ be a $3$-Kenmotsu manifold of dimension $\geq9$. Then, it is an Einstein space.
\end{Theorem}
\proof
	Theorem 3.3 in \cite{Ishihara74} implies that the quaternion Kahler manifold $\bar{M}$ is an Einstein space. Then, it follows from (\ref{ric}) that the Ricci tensor $Ric$ of $3$-Kenmotsu manifold $M$ satisfies the Einstein equation. Thus, $M$ is an Einstein space.\qed
Now, we define the $\varphi_\alpha-$holomorphic sectional curvature of $X\in\Gamma H$ by:
\begin{equation*}
\begin{split}
\mathbf{H}_\alpha(X)=& K(X,\varphi_\alpha X)\\
=&-\frac{R(X,\varphi_\alpha X,X,\varphi_\alpha X)}{g(X,X)g(\varphi_\alpha X,\varphi_\alpha X)-g(X,\varphi_\alpha X)^2}\\
=&-\frac{R(X,\varphi_\alpha X,X,\varphi_\alpha X)}{g(X,X)^2},
\end{split}
\end{equation*}
where $\alpha=1,2,3$.
\begin{Theorem}\label{H curv}
	The $\varphi_\alpha-$holomorphic sectional curvatures $\mathbf{H}_\alpha$ satisfy the following equation on a $3$-Kenmotsu manifold,
	$$\mathbf{H}_1(X)+\mathbf{H}_2(X)+\mathbf{H}_3(X)=-3,$$
	where $X\in\Gamma H$.
\end{Theorem}
\proof
	Let $(\varphi,\eta,\xi,g)$ be a Kenmotsu structure. Then, it follows from (\ref{1}), (\ref{2}) and Proposition \ref{Levi-Civita} that,
	$$R(X,Y)\varphi Z=g(\varphi Y,Z)X-g(\varphi X,Z)Y+g(Y,Z)\varphi X-g(X,Z)\varphi Y+\varphi(R(X,Y)Z).$$
	Therefore,
	\begin{equation*}
	\begin{split}
	R(X,Y,\varphi Z,\varphi W)=& g(R(X,Y)\varphi Z,\varphi W)=\\
	& g(\varphi Y,Z)g(X,\varphi W)-g(\varphi X,Z)g(Y,\varphi W)\\ &+g(Y,Z)g(\varphi X,\varphi W)-g(X,Z)g(\varphi Y,\varphi W)\\
	&+g(\varphi(R(X,Y)Z),\varphi W).
	\end{split}
	\end{equation*}
	Considering this equation on a $3$-Kenmotsu structure for each $\varphi_\alpha$. Let $\varphi=\varphi_1$, $Z=X$ and $Y=W=\varphi_3X$, then
	\begin{equation*}
	\begin{split}
	-R(X,\varphi_3X,\varphi_1X,\varphi_2X)=&-g(X,X)g(\varphi_2X,\varphi_2X)+R(X,\varphi_3X,X,\varphi_3X)\\
	=&-g(X,X)^2+R(X,\varphi_3X,X,\varphi_3X).
	\end{split}
	\end{equation*}
	Dividing both sides by $g(X,X)^2$, it implies:
	$$\frac{R(X,\varphi_3X,\varphi_2X,\varphi_1X)}{(g(X,X))^2}=-1-\mathbf{H}_3(X).$$
	Consider even permutations of $(1,2,3)$, one can get two other similar equations for $\mathbf{H}_1(X)$ and $\mathbf{H}_2(X)$. Then, the proof is completed by adding these three equations and using the Bianchi identity.\qed
%%%%%%%%%%%%%%%%%%%%%%%%%%%%%%%%%%%%%%%%%%%%%%%%%%%%%%%%%%%%%%%%%%%%%%%%%
\section{An Example of $3$-Kenmotsu Manifolds}
Consider the manifold $M$ given by $\{(x_0,x_1,x_2,x_3,x_4)\in I\!\!R^5\ |\ x_0\neq0\}$ and vector fields
$$\xi=X_0=-x_0\frac{\partial}{\partial x_0},\hspace{0.3cm}X_1=x_0\frac{\partial}{\partial x_1},\hspace{0.3cm}X_2=x_0\frac{\partial}{\partial x_2},\hspace{0.3cm}X_3=x_0\frac{\partial}{\partial x_3},\hspace{0.3cm}X_4=x_0\frac{\partial}{\partial x_4},$$
on $M$. It is easy to check that these vector fields satisfy following properties
$$[X_i,X_j]=0,\hspace{1cm}[\xi,X_i]=-X_i,$$
where $i,j=1,2,3,4$. Let $g$ be a Riemannian metric on $M$ such that,
$$g(X_i,X_j)=\delta_{ij}=\left\{\begin{array}{l}1\hspace{1cm}i=j\\0\hspace{1cm}i\neq j\end{array}\hspace{1cm}i,j=0,1,2,3,4\right.$$
Let $\varphi_\alpha$ for $\alpha=1,2,3$ be $(1,1)$-tensor field on $M$ defined by,
$$\begin{array}{l}
\varphi_1(\xi)=0,\hspace{0.3cm}\varphi_1(X_1)=X_2,\hspace{0.3cm}\varphi_1(X_2)=-X_1,\hspace{0.3cm}\varphi_1(X_3)=X_4,\hspace{0.3cm}\varphi_1(X_4)=-X_3,\\
\varphi_2(\xi)=0,\hspace{0.3cm}\varphi_2(X_1)=X_3,\hspace{0.3cm}\varphi_2(X_2)=-X_4,\hspace{0.3cm}\varphi_2(X_3)=-X_1,\hspace{0.3cm}\varphi_2(X_4)=X_2,\\
\varphi_3(\xi)=0,\hspace{0.3cm}\varphi_3(X_1)=X_4,\hspace{0.3cm}\varphi_3(X_2)=X_3,\hspace{0.3cm}\varphi_3(X_3)=-X_2,\hspace{0.3cm}\varphi_3(X_4)=-X_1.
\end{array}$$
It is easy to check for all $(i,j,k)$ as an even permutation of $(1,2,3)$ we have, $$\varphi_k=\varphi_i\circ\varphi_j.$$
If we define the $1$-form $\eta$ by $\eta(.)=g(\xi,.)$, it will be a straightforward calculation to check other properties in (\ref{0}). In this case, the corresponding Levi-Civita connection of $g$ on $M$ will have the following components,
\begin{equation}\label{ex levi}
\nabla_\xi\xi=\nabla_\xi X_i=\nabla_{X_i}\xi-X_i=\nabla_{X_i}X_j+\delta_{ij}\xi=0,
\end{equation}
where $i,j=1,2,3,4$. To show $(\varphi_\alpha,\eta,\xi,g)$, for $\alpha=1,2,3$, is a $3$-Kenmotsu structure one should check equation (\ref{1}) which is easy by the help of (\ref{ex levi}) and definition of $\varphi_\alpha$ where $\alpha=1,2,3$.\par
Moreover, (\ref{ex levi}) implies that components of the curvature tensor $R$ satisfy the following equations:
\begin{alignat*}{2}
& R(X_i,X_j)\xi=0 && i,j=1,2,3,4,\\
& R(\xi,X_i)X_j=-R(X_i,\xi)X_j=\nabla_{X_i}X_j && i,j=0,1,2,3,4,\\
& R(X_i,X_j)X_k=g(X_i,X_k)X_j-g(X_j,X_k)X_i\ \ \ \ \ && i,j=1,2,3,4.
\end{alignat*}
Since $\{X_1,X_2,X_3,X_4\}$ is an orthonormal basis of sections on $H$, we can write all $X\in\Gamma H$ as $X=\sum_{i=1}^4a_iX_i$, where $a_i$ is a scalar function on $M$ for $i=1,2,3,4$. Then, the $\varphi_\alpha-$holomorphic sectional curvature $\mathbf{H}_\alpha(X)$ of $X\in\Gamma H$ for $\alpha=1,2,3$ can be calculated as follows:
\begin{align*}
\mathbf{H}_\alpha(X)=&-\frac{R(X,\varphi_\alpha X,X,\varphi_\alpha X)}{g(X,X)g(\varphi_\alpha X,\varphi_\alpha X)-g(X,\varphi_\alpha X)^2}\\
=&-\frac{g(R(X,\varphi_\alpha X)X,\varphi_\alpha X)}{g(X,X)^2-g(X,\varphi_\alpha X)^2}\\
=&-\frac{g(g(X,X)\varphi_\alpha X-g(\varphi_\alpha X,X)X,\varphi_\alpha X)}{g(X,X)^2-g(X,\varphi_\alpha X)^2}=-1
\end{align*}
which confirms the result of Theorem \ref{H curv}.
%%%%%%%%%%%%%%%%%%%%%%%%%%%%%%%%%%%%%%%%%%%%%%%%%%%%%%%%%%%%%%%%%%%%%%%%%

\end{document}